\documentclass[12pt]{amsart}
\usepackage{amssymb}

\newtheorem{theorem}{Theorem}

\newtheorem{proposition}{Proposition}
\newtheorem{corollary}{Corollary}
\newtheorem{remark}{Remark}

\newtheorem{lemma}{Lemma}

\newcommand{\Z}{{\mathbb Z}}
\newcommand{\Q}{{\mathbb Q}}
\newcommand{\R}{{\mathbb R}}
\newcommand{\C}{{\mathbb C}}

\newcommand{\RP}{{\mathbb RP}}

\begin{document}

\title[]{Homology of the complex of not 2-divisible partitions}

\author{V.A.~Vassiliev}

\date{}

\address{Steklov Mathematical Institute and \newline 
National Research University Higher School of Economics, Moscow}

\email{vva@mi.ras.ru}

\thanks{}

\dedicatory{}



\begin{abstract}
The rational homology group of the order complex of non-even partitions of a finite set is calculated. A twisted version of the Goresky-MacPherson approach to similar homology calculations is proposed.
\end{abstract}

\maketitle

\section{Introduction}
The partitions of a $N$-element set form a partially ordered set (any partition dominates its subpartitions); so its order complex is well-defined. Its subcomplex of {\it 2-divisible partitions} (into the parts of even cardinalities) was explored in \cite{Syl} in a connection with the statistical physics; for some generalizations see \cite{Sta}, \cite{ER}. 
We calculate the rational homology group of the order complex of the complementary set of non-complete partitions having several parts of odd cardinality.
\medskip

{\bf Notation.} 
For any natural $N\geq 3$, denote by $\Delta(N)$ the order complex of the poset of all partitions of the set $\{1, \dots, N\}$, having at least one part of cardinality $\ge 2$. $\partial \Delta(N)$ is the {\em link} of $\Delta(N)$, that is, the subcomplex consisting of all simplices not containing its maximal vertex (corresponding to the ``partition'' into the single set).
For any partition $A$, $\Delta_A \subset \Delta(N)$ is the order complex of the poset of subpartitions of $A$ (including $A$ itself); 
$\partial \Delta_A$ is the link of $\Delta_A$. 

If $N$ is even, denote by $\Xi_2(N) \subset \Delta(N)$ the order complex of the poset of partitions of $\{1, \dots, N\}$ into parts, some of which consist of odd numbers of elements.

$\tilde H_*$ denotes the homology groups reduced modulo a point, and $\bar H_*$ the homology groups of complexes of locally finite singular chains.
\medskip

\begin{theorem}
\label{main}
For any even $N>2$, the group $H_i(\Xi_2(N),\Q)$ is trivial for all positive $i \neq N-3$, and the dimension of $H_{N-3}(\Xi_2(N),\Q)$ is equal to $(N-1)!! (N-3)!!$ .
\end{theorem}

\begin{remark} \rm
This problem comes from the equivariant homotopy theory: the groups in question arise in the calculation of the homology groups of spaces of odd maps of spheres, see \cite{dan18}.
\end{remark} 
 
There is a method of computing the homology of order complexes, based on the Goresky-MacPherson formula \cite{GM}.
Originally, this formula expresses the cohomology groups of complements of plane arrangements through the homology of order complexes of these planes and their intersections. However, the cohomology of complements of arrangements often can be calculated also by easier techniques. This allows one to use this formula in the opposite direction: see \S \ref{triv} for its application to a standard example.

We propose a similar method exploiting the homology groups with coefficients in non-trivial local systems. I hope that it will work for other problems of this kind, e.g. for the complexes of not $k$-divisible partitions (which should be related with the ``$k$-equal arrangements'').

\section{A preparatory example: the complex of non-complete partitions}
\label{triv}

\begin{proposition}[see e.g. \cite{BW}, \cite{StaBook}]
\label{trii}
$\tilde H_i(\partial \Delta(N)) =0$ if $i \neq N-3$, and $$\tilde H_{N-3} (\partial \Delta(N)) \simeq \Z^{(N-1)!}.$$
\end{proposition}

\noindent
{\it Proof} (cf. \cite{conn}). Let $I(\R^m,N)$ be the space of ordered configurations of $N$ distinct points in $\R^m$. This is an open subset in $\R^{mN}$; its complement is a subspace arrangement whose strata are in the one-to-one correspondence with the partitions of $\{1, \dots, N\}$ into $\leq N-1$ subsets. The Goresky-MacPherson formula says that 
\begin{equation}
\label{GMf}
\tilde H^i(I(\R^m,N)) \simeq \bigoplus_{A} \tilde H_{mN - m k(A) - i-2} (\partial \Delta_A),
\end{equation} 
where the summation is taken over all such partitions $A$ of $\{1, \dots, N\}$ and $k(A)$ is the number of parts of the partition $A$ (the group $\tilde H_*(\partial \Delta_A)$ is assumed to be equal to $\Z$ in dimension $-1$ only if $\Delta_A$ is a point). 

 On the other hand, the cohomology group of $I(\R^m,N)$ can be easily calculated as in \cite{fuchs}: it has Poincar\'e polynomial $\prod_{j=1}^{N-1}(1+jt^{m-1})$ and no torsion. In particular, the group $H^i(I(\R^m,N))$ is equal to $\Z^{(N-1)!}$ for $i=(N-1)(m-1)$ and is trivial if $i>(N-1)(m-1)$. 

If $k(A)\geq 2$ then the right-hand group in (\ref{GMf}) is trivial for $i > mN-2m -2$. Therefore if $ m>N-3$ (and hence $(N-1)(m-1)> mN-2m-2)$, then only the homology groups of $\partial \Delta_A$ corresponding to the ``partition'' $A$ into the single set $\{1, \dots, N\}$ can contribute by (\ref{GMf}) to the groups $\tilde H^i(I(\R^m,N))$ with $i \geq (N-1)(m-1)$. Namely, we get $\tilde H_{N-3}(\partial \Delta_A) \simeq \tilde H^{(m-1)(N-1)}(I(\R^m,N)) \simeq \Z^{(N-1)!}$, and $\tilde H_j(\partial \Delta_A)=0$ for all $j<N-3$. Since $\Delta_A \equiv \Delta(N)$ for this ``partition'', Proposition is proved.
\hfill $\Box$

\section{The twisted modification}

Now we assume again that $N$ is even. Consider the projective space $\RP^m$ with odd $m>1$, and the local system of groups on it, locally isomorphic to the constant sheaf with fiber $\Q$, but with the non-trivial monodromy action of $\pi_1(\RP^1) \sim \Z_2$: the unique non-zero element of this group acts on the fibers as multiplication by $-1$. The homology group of $\RP^m$ with coefficients in this system is trivial in all dimensions including 0. Consider also the space $(\RP^m)^{N-1}$, i.e.~the Cartesian power of $N-1$ copies of $\RP^m$. Denote by $\Theta$ the local system on it which is the tensor product of $N-1$ systems lifted (by the standard projections) from the above-considered systems on its factors. 

Distinguish a point $x_1 \in \RP^m$. Denote by $\Phi(\RP^m,N)$ the subset in $(\RP^m)^{N-1}$ consisting of points $(x_2, \dots, x_{N})$, $x_i \in \RP^m$, such that all these points are different from one another and from $x_1$. Denote by $\Sigma$ 
its complement $(\RP^m)^{N-1}\setminus \Phi(\RP^m,N)$. Consider the group 
\begin{equation}
\label{abs}
H^*(\Phi(\RP^m,N), \Theta)
\end{equation}
 and its Poincar\'e--Lefschetz dual group 
\begin{equation}
\label{relat}
\bar H_*(\Phi(\RP^m,N), \Theta) \simeq H_*((\RP^m)^{N-1}, \Sigma; \Theta)\ .
\end{equation} 
 Since $H_*((\RP^m)^{N-1}, \Theta) \equiv 0$, the boundary map defines an isomorphism of the group (\ref{relat}) to $H_*(\Sigma, \Theta)$. 

\begin{proposition}
\label{lem1}
The group $($\ref{abs}$)$ is trivial in all dimensions not equal to \ $a(m-1)$ for some $a\in \{1,2, \dots, N-1\}.$
\end{proposition}

{\it Proof.} The universal covering of $\Phi(\RP^m,N)$ can be described as the space of all sequences $(\tilde x_2, \dots, \tilde x_{N}) \subset S^m$ such that $\tilde x_2$ does not coincide with either of two preimages of the point $x_1$ under the obvious projection $S^m \to \RP^m$, $\tilde x_3$ does not coincide with either of four points: these preimages of $x_1$, the point $\tilde x_2$, and its antipode; $\tilde x_4$ does not coincide with either of five points defined on the previous steps and with the antipode of $\tilde x_3,$ etc. The space of all such sequences obviously is a tower of fiber bundles,
all whose fibers are homotopy equivalent to wedges of $(m-1)$-dimensional spheres. By the standard spectral sequence considerations, its homology group may be non-trivial only in the dimensions $0, (m-1), 2(m-1), \dots, (N-1)(m-1)$. 

Consider the following action of the group $(\Z_2)^{N-1}$ on this covering space: the non-trivial element of the $i$th summand of this group sends $\tilde x_{i+1}$ to its antipode and leaves all other points $\tilde x_j$ unmoved. Consider also the action of this group $(\Z_2)^{N-1}$ on the group of singular chains of this covering space: any element of $(\Z_2)^{N-1}$ moves any simplex geometrically in accordance with the previous action, and additionally multiplies it by $\pm 1$ depending on the parity of the number of non-trivial components in this element. The complex of chains of $\Phi(\RP^m,N)$ with coefficients in $\Theta$ is the quotient of the complex of ordinary chains (with rational coefficients) of the covering space by this action. Therefore by Proposition 2 from Appendix to \cite{serre} the group (\ref{abs}) is a quotient group of the homology group of this covering space. In addition, it is obviously trivial in dimension 0. \hfill $\Box$

\begin{corollary}
\label{corr2}
The groups $H_j(\Sigma, \Theta)$ can be non-trivial only for $j$ of the form $N-2 + a(m-1),$ $a \in \{0, 1, \dots, N-2\}.$ 
\end{corollary}

\noindent
{\it Proof:} the boundary and Poincar\'e--Lefschetz isomorphisms $$H_j(\Sigma,\Theta) \gets H_{j+1}((\RP^m)^{N-1},\Sigma; \Theta) \simeq H^{m(N-1)-(j+1)}(\Phi (\RP^m,N), \Theta). \hfill \Box$$

\section{Simplicial resolution and spectral sequence}

The subset $\Sigma \subset (\RP^m)^{N-1}$ has the following standard simplicial resolution $\sigma \subset \Delta(N) \times (\RP^m)^{N-1} .$ For any partition $A$ of the set $\{1, \dots, N\}$ into $k$ parts $A_1, \dots, A_k$ we consider the subspace $L(A) \subset (\RP^m)^{N-1}$ consisting of all points $(x_2, \dots, x_N)$ such that all points $x_i\in \RP^m$ (including $x_1$), whose indices $i$ belong to one and the same part $A_j$, do coincide. Obviously, in this case $L(A) \sim (\RP^m)^{k-1}$. 

The space $\sigma$ is defined as the union of all spaces $\Delta_A \times L(A)$ over all partitions $A$ of $\{1, \dots, N\}$ into $\leq N-1$ parts.

Let $\pi: \sigma \to \Sigma$ be the map induced by the standard projection $\Delta(N) \times (\RP^m)^{N-1} \to (\RP^m)^{N-1}$, and $\tilde \Theta$ the local system on $\sigma$ induced from $\Theta$ by this map.

\begin{proposition}
\label{leray}
 The map $\pi$ induces an isomorphism
$H_*(\sigma, \tilde \Theta) \simeq H_*(\Sigma, \Theta).$
\end{proposition}

\noindent 
{\it Proof.}
For any $z \in \Sigma$, the fiber $\pi^{-1}(z)$ is contractible and is a deformation retract of the preimage $\pi^{-1}(U)$ of a neighborhood $U$ of $z$. Therefore the direct images $R^q_\pi(\tilde \Theta)$ of the sheaf $\tilde \Theta$ are trivial for $q>0$, and the standard morphism $R^0_\pi(\tilde \Theta) \to \Theta$ is an isomorphism. The  Proposition  follows by the Leray spectral sequence of $\pi$, see \S 4.17 in \cite{Godem}. \hfill $\Box$ \medskip

The space $\sigma$ has a natural filtration $\sigma_0 \subset \sigma_1 \subset \dots \subset \sigma_{N-2} = \sigma:$ the subset $\sigma_p$ is the union of all products $\Delta_A \times L(A)$ over the partitions $A$ into $\geq N-p-1$ parts. Let $\{E_{p,q}^r\}$ be the spectral sequence calculating the group $H_*(\sigma, \tilde \Theta)$ and generated by this filtration. Then $E^1_{p,q} \simeq \bar H_{p+q} (\sigma_p \setminus \sigma_{p-1}, \tilde \Theta)$.

For any partition $A$, denote by $\breve \Delta_A$ the space $\Delta_A \setminus \partial \Delta_A$, so that $H_i(\partial \Delta_A) \simeq \bar H_{i+1}(\breve \Delta_A) \simeq H_{i+1}(\Delta_A, \partial \Delta_A)$ for any $i$.
The set $\sigma_p \setminus \sigma_{p-1}$ is the disjoint union of spaces corresponding to all partitions of $\{1, \dots, N\}$ into exactly $N-p-1$ parts. Any such space corresponding to some partition $A$ is homeomorphic to $\breve \Delta_A \times L(A)$, and the restriction of the sheaf $\tilde \Theta$ to it is equal to that lifted from the sheaf $\Theta$ on $L(A)$.

\begin{proposition}
\label{proo}
For any partition $A$ of the set $\{1, \dots, N\}$ into $k<N$ parts $A_1, \dots, A_k$, the group $\bar H_*(\breve \Delta_A \times L(A), \tilde \Theta)$ is 
trivial if $A$ is not 2-divisible and equals $\bar H_*(\breve \Delta_A,\Q) \otimes H_*((\RP^m)^{k-1}, \Q)$ if $A$ is 2-divisible. 
\end{proposition}

\noindent
{\it Proof.} If the partition $A$ is 2-divisible, then the restriction of $\Theta$ to $L(A)$ is the constant $\Q$-sheaf, and $H_*(L(A),\Theta) \simeq H_*((\RP^m)^{k-1},\Q).$ Otherwise it is a tensor product of several local systems on the factors $\RP^m$, some of which are non-trivial, hence $H_*(L(A),\Theta) \equiv 0$. The K\"unneth formula completes the proof. \hfill $\Box$ \medskip

Let $\Xi_2(\RP^m,N) \subset \sigma$ be the union of all blocks $\breve \Delta_A \times L(A)$ over the not 2-divisible partitions $A$. It is a compact subset in $\sigma$.

\begin{corollary}
\label{propop8}
$H_*(\Xi_2(\RP^m,N), \tilde \Theta) \equiv 0.$ 
The inclusion \ $\sigma \setminus \Xi_2(\RP^m,N) \hookrightarrow \sigma$ induces an isomorphism of homology groups with closed supports, \smallskip

\hspace{3cm} $\bar H_*(\sigma \setminus \Xi_2(\RP^m,N),\tilde \Theta) \simeq H_*(\sigma,\tilde \Theta).$ \hfill $\Box$
\end{corollary} 

\begin{corollary}
\label{corr5}
The groups $E^1_{p,q}$ of our spectral sequence can be non-trivial only for the pairs $(p,q)$ satisfying the following conditions:

1. $p \in \left[\frac{N}{2}-1, N-2\right],$

2. $q \in \{0, m, 2m, \dots , (\frac{N}{2}-1)m\}$,

3. $ pm+q \leq (N-2)m .$ \hfill $\Box$
\end{corollary}

\unitlength=0,6mm
\begin{figure}
$$
\mbox{
\begin{picture}(110,87)
\put(10,7){\vector(0,1){80}}
\put(10,7){\vector(1,0){60}}
\put(20,7){\line(0,1){70}}
\put(30,7){\line(0,1){70}}
\put(40,7){\line(0,1){70}}
\put(50,7){\line(0,1){70}}
\put(60,7){\line(0,1){70}} 
\put(10,17){\line(1,0){50}}
\put(10,27){\line(1,0){50}}
\put(10,37){\line(1,0){50}}
\put(10,47){\line(1,0){50}}
\put(10,57){\line(1,0){50}}
\put(10,67){\line(1,0){50}}
\put(10,77){\line(1,0){50}}
\put(4,10){$0$}
\put(3,40){$m$}
\put(0,70){$2m$}
\put(5,82){$q$}
\put(15,0){$0$}
\put(25,0){$1$}
\put(35,0){$2$}
\put(45,0){$3$}
\put(55,0){$4$}
\put(67,2){$p$}
\put(32,10){15}
\put(42,10){90}
\put(52,10){$5!$}
\put(32,40){30}
\put(42,40){90}
\put(32,70){15}
\end{picture}
}
\mbox{
\begin{picture}(70,87)
\put(10,7){\vector(0,1){80}}
\put(10,7){\vector(1,0){60}}
\put(20,7){\line(0,1){70}}
\put(30,7){\line(0,1){70}}
\put(40,7){\line(0,1){70}}
\put(50,7){\line(0,1){70}}
\put(60,7){\line(0,1){70}} 
\put(10,17){\line(1,0){50}}
\put(10,27){\line(1,0){50}}
\put(10,37){\line(1,0){50}}
\put(10,47){\line(1,0){50}}
\put(10,57){\line(1,0){50}}
\put(10,67){\line(1,0){50}}
\put(10,77){\line(1,0){50}}
\put(4,10){$0$}
\put(3,40){$m$}
\put(0,70){$2m$}
\put(5,82){$q$}
\put(15,0){$0$}
\put(25,0){$1$}
\put(35,0){$2$}
\put(45,0){$3$}
\put(55,0){$4$}
\put(67,2){$p$}
\put(52,10){45}
\put(42,40){60}
\put(32,70){15}
\end{picture}
}
$$
\caption{Pages $E^1$ and $E^2$ for $N=6$}
\label{fig1}
\end{figure}
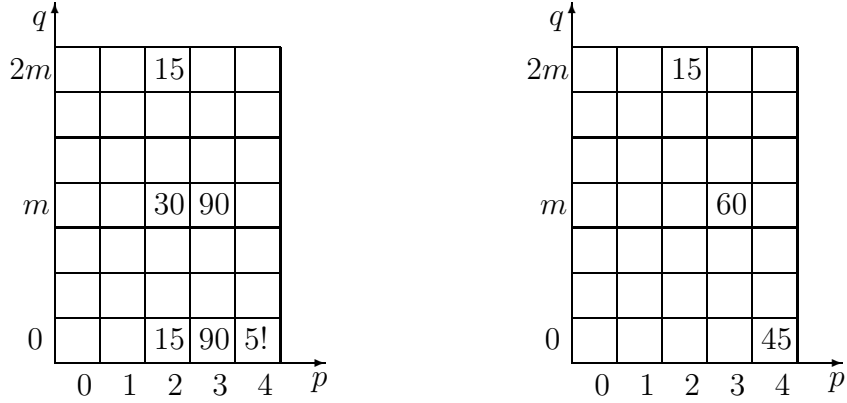

So, these non-trivial groups $E^1_{p,q}$ lie in the triangle with the vertices $\left(\frac{N}{2}-1, 0\right), (N-2,0), \left(\frac{N}{2}-1, (\frac{N}{2}-1)m\right).$ The term $(E^1,d^1)$ of the spectral sequence consists of $\frac{N}{2}$ horizontal complexes placed on the lines $q=0, m, \dots, (\frac{N}{2}-1)m$ and consisting of $\frac{N}{2}, \frac{N}{2}-1, \dots, 1$ non-trivial groups respectively. These groups for $m=3$ and $N=6$ are shown in Fig.~\ref{fig1} left (where only the ranks of the groups are indicated, and all unfilled cells mean the trivial groups).

\begin{lemma}
\label{leme2}
If $m > \frac{N}{2}-1$ then 
all these $\frac{N}{2}$ horizontal complexes are acyclic in all their groups except for the top-dimensional ones $($i.e., only the groups $E^2_{N-2,0},$ $E^2_{N-3, m},$ $E^2_{N-4, 2m}, \dots , E^2_{\frac{N}{2}-1,(\frac{N}{2}-1)m}$ can be non-trivial, see Fig.~\ref{fig1} right$)$.
\end{lemma} 

\noindent
{\it Proof.} By Corollary \ref{corr5} and dimensional reasons, $E^2 \equiv E^\infty$ for such $m$, and all non-trivial groups $E^2_{p,q}$ are equal to  $H_{p+q}(\sigma, \tilde \Theta) \simeq H_{p+q}(\Sigma, \Theta)$. 
But only those groups of our complexes, which are mentioned in Lemma, have the numbers $p+q$ as indicated in Corollary \ref{corr2}.\hfill $\Box$ \smallskip

(In fact, the assertion of this lemma is true for arbitrary odd $m$, but the proof is then not so immediate). \medskip

The order complex $\Delta(N)$ of partitions of $\{1, \dots, N\}$ is embedded into $\sigma$ as the preimage of the point $(x_1, x_1, \dots, x_1) \in \Sigma$ under the projection $\sigma \to \Sigma$: indeed, the restriction of the projection $\Delta(N) \times (\RP^m)^{N-1} \to \Delta(N)$ to this preimage is a homeomorphism. The subcomplex $\Xi_2(N)$ of $\Delta(N)$ is the preimage of $\Xi_2(\RP^m,N)$ under this embedding. 

\begin{proposition}
\label{propo9}
If $m> \frac{N}{2}-1$ then this embedding induces the isomorphisms $\bar H_j(\Delta(N) \setminus \Xi_2(N),\Q) \simeq \bar H_j(\sigma \setminus \Xi_2(\RP^m,N), \tilde \Theta)$ in all dimensions $j < N-3+m$. 
\end{proposition}

\noindent
{\it Proof.} The preimage of the term $\sigma_p$ of our filtration on $\sigma$ under this embedding is exactly the union of subspaces $\breve \Delta_A \subset \Delta(N)$ over the partitions into $\geq N-p-1$ parts. Let $\{{\mathcal E}^r_{p,q}\}$ be the homological spectral sequence calculating $\bar H_*(\Delta(N) \setminus \Xi_2(N),\Q)$ and induced by this filtration. By Proposition \ref{trii} its terms ${\mathcal E}^1_{p,q}$ can be non-trivial for $q=0$ only. By Proposition \ref{proo}, the induced morphism ${\mathcal E}^1 \to E^1$ of first terms of our spectral sequences is an isomorphism onto the bottom line $\{q=0\}$. Therefore and by Lemma \ref{leme2}, the induced morphisms ${\mathcal E}^\infty_{p,q} \to E^\infty_{p,q}$ of limit groups are isomorphisms for all $p,q$ with $p+q <N-3+m$.  \hfill $\Box$

\begin{corollary}
$H_j(\Xi_2(N), \Q) = 0$ for $j <N-3$. 
\end{corollary}

\noindent
{\it Proof.} We have $H_j(\Xi_2(N), \Q) \simeq \bar H_{j+1}(\Delta(N) \setminus \Xi_2(N),\Q) \simeq \bar H_{j+1}(\sigma \setminus \Xi_2(\RP^m,N), \tilde \Theta) \simeq H_{j+1}(\sigma, \tilde \Theta) \simeq H_{j+1}(\Sigma, \Theta) \simeq 0$
for such $j$ and sufficiently large $m$: here the first equality is the boundary isomorphism in $\Delta(N)$, and the next four ones follow respectively from  
Proposition \ref{propo9}, Corollary \ref{propop8}, Proposition \ref{leray}, and Corollary \ref{corr2}.\hfill $\Box$ 

\section{Euler characteristic}
\label{euler}

To prove Theorem \ref{main}, it remains to check that the Euler characteristic of the subcomplex $\Xi_2(N) \subset \Delta(N)$ is equal to $-(N-1)!!(N-3)!!+1$ or, which is the same, the Euler characteristic (with closed supports) of the difference $\Delta(N) \setminus \Xi_2(N)$ is equal to $(N-1)!!(N-3)!!$. The latter space is represented by all simplices of $\Delta(N)$, whose maximal vertices correspond to 2-divisible partitions of $\{1, \dots, N\}$. Given such a partition, let $A_1, \dots, A_k$ be its parts ordered by the values of their smallest elements, and $a_1, \dots, a_k$ the cardinalities of these parts. By Proposition \ref{trii}, the union of simplices, whose maximal vertex corresponds to this partition, contributes the number $(-1)^{k+1}\prod_{j=1}^k (a_j-1)!$ to the desired Euler characteristic. The absolute value of this number can be interpreted as the number of reorderings $(\alpha_1, \dots, \alpha_N)$ of $\{1, \dots, N\}$ having the following form: the first element $\alpha_1$ is the smallest element of the part $A_1$ (it can be the element $\{1\}$ only), then all other elements of $A_1$ follow in an arbitrary order, then we take  the smallest element of $A_2$  for $\alpha_{a_1+1}$, then all other elements of $A_2$ ordered arbitrarily, etc. 

Now consider an arbitrary permutation $(\alpha_1, \dots, \alpha_N)$ of $\{1, \dots, N\}$ starting with $\{1\}$, and count the occurencies of this permutation from the even partitions by the above construction. If for no odd number $r>1$ the corresponding element $\alpha_r$ is smaller than all the remaining elements $\alpha_{r+1}, \dots, \alpha_N$, then this permutation occurs exactly once from the maximal ``partition'' into only one set. If however for some $t$ different numbers $r_s \in \{3, 5, \dots, N-1\}$ the elements $\alpha_{r_s}$ are smaller than all the forthcoming ones, then this permutation occurs once from the maximal partition, $t$ times from the partitions into two even parts, $\binom{t}{2}$ times from partitions into $3$ parts, ..., $t$ times from the partitions into $t$ parts, and once from a partition into $t+1$ parts. These occurencies should be counted with alternating signs in the calculation of the Euler characteristic, and hence annihilate one another. So, only the permutations without such minimal elements $\alpha_r$ should be accounted. The number of them is equal to $(N-1)!!(N-3)!!$ because enumerating them we have one value prohibited on any even step. \hfill $\Box$

\section{Concluding questions and remarks}

1. Theorem \ref{main} can be extended almost immediately to the homology with coefficients in any ring with division by 2, but is it true also over the integers? 

2. Find a more conventional proof for it (say, by the method of \cite{forman}).
Maybe the way of combining different cycles described in \S \ref{euler} can give a hint for the reduction of simplices participating in such a proof ?

3. It looks likely that the above approach can provide a useful heuristics for other problems of this kind, for instance for the structures related with the ``$k$-equal manifolds'', see \cite{BW}, and the $k$-divisible partitions of \cite{Sta}.
Probably it is useful in this case to replace the projective space $\RP^m$ by the lens space $S^{2n-1} / \Z_{k}$ and the local system $\Theta$ by a similar local system with the fiber $\C$ and the $\pi_1$-action defined by the multiplication by the powers of $e^{2\pi i/k}$.

4. The entire group (\ref{abs}) is also calculated: its Poincar\'e polynomial is equal to 
$$(N-1)!! \ t^{\frac{N}{2}(m-1)} \prod_{j=1}^{\frac{N}{2}-1} \left(1+(2j-1)t^{m-1}\right).$$

\end{document}